\documentclass[12pt]{article}
\usepackage[utf8]{inputenc}
\usepackage{amssymb,amsmath,amsfonts, eucal, mathrsfs, amsthm} 
\usepackage[colorinlistoftodos]{todonotes}
\usepackage[allcolors=blue]{hyperref}
\newtheorem{theorem}{Theorem}
\newtheorem{proposition}[theorem]{Proposition}
\newtheorem{lemma}[theorem]{Lemma}

\newtheorem{corollary}[theorem]{Corollary}

\newtheorem{remark}[theorem]{Remark}
\newtheorem{remarks}[theorem]{Remarks}

\newtheorem{examples}[theorem]{Examples}

\theoremstyle{definition}
\newcommand{\R}{\mathbb{R}}
\newcommand{\Q}{\mathbb{Q}}
\newcommand{\Sf}{\mathbb{S}}

\newcommand{\Hy}{\mathbb{H}}
\newcommand{\spa}{\mbox{span}}

\newcommand{\Ima}{\mbox{Im }}

\newcommand{\nap}{\nabla^{\perp}}
\newcommand{\nab}{\tilde\nabla}

\newcommand{\Hom}{\mbox{Hom}}

\newcommand{\Y}{\mathcal{Y}\,}
\def\<{{\langle}}
\def\>{{\rangle}}

\def\F{{\cal F}}

\def\Sal{{\cal S}}

\def\Y{{\cal Y}}

\def\d{\partial}
\def\a{\alpha}

\def\id{I}

\def\be{\begin{equation} }
\def\ee{\end{equation} }
\def\Ima{Im}
\begin{document}

\date{}
\title{Genuine infinitesimal bendings\\ of submanifolds}
\author{M. Dajczer and M. I. Jimenez}
\maketitle

\begin{abstract} A basic question in submanifold theory is 
whether a given isometric immersion $f\colon M^n\to\R^{n+p}$ of 
a Riemannian manifold of dimension $n\geq 3$ into Euclidean space 
with low codimension $p$ admits, locally or globally, a genuine 
infinitesimal bending. That is, if there exists a genuine smooth 
variation of $f$  by immersions that are isometric up to the first 
order. Until now only the hypersurface case $p=1$ was well understood. 
We show that a strong necessary local condition to admit such a bending 
is the submanifold to be ruled and give a lower bound for the dimension 
of the rulings. In the global case, we describe the situation of 
compact submanifolds of dimension $n\geq 5$ in codimension $p=2$.   
\end{abstract}

An isometric immersion $f\colon M^n\to\R^{n+p}$ of an $n$-dimensional 
Riemannian manifold $M^n$ into Euclidean space with codimension $p$ is 
called isometrically bendable if there is a non-trivial smooth variation 
$\F\colon I\times M^n\to\R^{n+p}$ of $f$ for an interval $0\in I\subset\R$ 
such that $f_t=\F(t,\cdot)\colon M^n\to\R^{n+p}$  with $f_0=f$ is an isometric 
immersion for any $t\in I$ , that is, the metrics $g_t$ induced by $f_t$ 
satisfy $g_t=g_0$. The bending being trivial means that the variation is 
the restriction to the submanifold of a smooth one-parameter family of 
isometries of $\R^{n+p}$.

The study of bendings of surfaces $M^2$ in $\R^3$ was a hot topic between 
geometers in the $19^{th}$ century.  Initially, there was no distinction 
between isometric variations and the ones that are only infinitesimally 
isometric, but that changed due to the work of Darboux by the end of 
that century. For a modern account of some aspects of the subject we 
refer to Spivak \cite{Sp}.  

The study of isometric bendings of hypersurfaces 
$f\colon M^n\to\R^{n+1}$, $n\geq 3$, goes back to the first part of 
the last century. In fact,  the local classification of isometrically 
bendable hypersurfaces is due to Sbrana \cite{Sb1} in 1909 and Cartan 
\cite{Ca} in 1916. For a modern presentation of their parametric 
classifications, as well as for further results, see \cite{DFT} or 
\cite{DT}. In the global case, the classification is due to
Sacksteder \cite{Sa} for compact hypersurfaces and to Dajczer and 
Gromoll \cite{DG0} in the case of complete hypersurfaces.

The classical concept of an infinitesimal bending of an isometric
immersion $f\colon M^n\to\R^{n+p}$ is the infinitesimal analogue 
of an isometric bending and refers to smooth variations 
$\F\colon I\times M^n\to\R^{n+p}$ that preserve lengths ``up to the 
first order", that is, the metrics $g_t$ induced by 
$f_t=\F(t,\cdot)\colon M^n\to\R^{n+p}$ satisfy $g'_t(0)=0$. 
The variational vector field $\tau=\F_*\d/\d t|_{t=0}$ verifies
\be\label{infbend}
\<f_*X,\tau_*X\>=0
\ee
for any tangent vector fields $X\in\mathfrak{X}(M)$. Clearly 
\eqref{infbend} is the condition for a smooth variation to preserve 
the metric up to the first order.  If $\tau$ is an immersion,
it was said classically that the pair of submanifolds 
$f$ and $\tau$ correspond with orthogonality of corresponding
linear elements; see Bianchi \cite{Bi} or Eisenhart \cite{Ei}. 

We say that a section $\tau$ of $f^*T\R^{n+p}$ is an \emph{infinitesimal 
bending} of an isometric immersion $f\colon M^n\to\R^{n+p}$ if 
\eqref{infbend} holds. Given a smooth variation whose variational vector 
field $\tau$ is an infinitesimal bending, by keeping only the terms 
of first order of the variation we obtain the smooth variation 
\mbox{$\F\colon\R\times M^n\to\R^{n+p}$} with variational vector 
field $\tau$ defined by $f_t=f+t\tau$. 
Then \eqref{infbend} gives
$$
\|f_{t*}X\|^2=\|f_*X\|^2+t^2\|\tau_*X\|^2
$$
for any $X\in TM$.

Of course, we always have the \emph{trivial infinitesimal 
bendings} obtained as the variational vector field of a smooth 
variation by isometries of the ambient space. In other words, 
they are locally the restriction to the submanifold of a Killing 
vector field of the ambient space.

Dajczer and Rodríguez \cite{DR} showed that  submanifolds  
in low codimension are generically \emph{infinitesimally rigid},
that is, only trivial  infinitesimal bendings are possible. 
In fact, they proved that well-known algebraic conditions on the 
second fundamental form of an immersion that give isometric rigidity 
also yield infinitesimal rigidity. For instance, for a hypersurface 
$f\colon M^n\to\R^{n+1}$ to be infinitesimally bendable it is a 
necessary condition (but far from sufficient) to have at most two 
nonzero principal curvatures at any point. This result is already 
contained in the book of Cesàro \cite{Ce} published in 1886. For 
higher codimension the rather strong algebraic conditions are given 
in terms of the type number or the $s$-nullities of the immersion.

After the pioneering work of Sbrana \cite{Sb} in 1908, a complete 
parametric local classification of the infinitesimally bendable 
hypersurfaces was given by Dajczer and Vlachos \cite{DV}.  
In particular, they showed that this class is much larger than 
the class of isometrically bendable ones, a fact that may be seen 
as a surprise. The classification in the case of complete 
hypersurfaces was obtained by Jimenez \cite{Ji}. 
Infinitesimal bendings of submanifolds have also been considered 
by Schouten \cite{sc} in 1928. 

When trying to understand the geometry of the infinitesimally 
bendable submanifolds in codimension larger than one the following 
fact has to be taken into consideration. If $\tilde{\tau}$ is an 
infinitesimal bending of an isometric immersion 
$F\colon\tilde{M}^{n+\ell}\to\R^{n+p}$, $0<\ell< p$, and 
$j\colon M^n\to\tilde{M}^{n+\ell}$ is an embedding, then 
$\tau=\tilde{\tau}|_{j(M)}$ is an infinitesimal bending of 
$f=F\circ j\colon M^n\to\R^{n+p}$.  This basic observation 
motivates the following definitions where a more general 
situation is considered since certain singularities are 
allowed.

A smooth map $F\colon\tilde{M}^{n+\ell}\to\R^{n+p}$, $0<\ell<p$, 
 from a differentiable manifold into Euclidean space is said 
to be a \emph{singular extension} of a given isometric immersion 
$f\colon M^n\to\R^{n+p}$ if there is an embedding 
\mbox{$j\colon M^n\to\tilde{M}^{n+\ell}$}, $0<\ell< p$, such that 
$F$ is an immersion along $\tilde{M}^{n+\ell}\setminus j(M)$ and 
$f=F\circ j$.
Notice that the map $F$ may fail (but not necessarily) to be an 
immersion along points of $j(M)$. We say that an infinitesimal 
bending $\tau$ of an isometric immersion $f\colon M^n\to\R^{n+p}$ 
\emph{extends in the singular sense} if there is a singular extension 
\mbox{$F\colon\tilde{M}^{n+\ell}\to\R^{n+p}$} of $f$ and a smooth 
map $\tilde{\tau}\colon\tilde{M}^{n+\ell}\to\R^{n+p}$ 
such that $\tilde{\tau}$ is an infinitesimal bending of 
$F_{\tilde{M}\setminus j(M)}$ and $\tau=\tilde{\tau}|_{j(M)}$.
\medskip

We point out that the necessity to admit the existence of 
singularities of $F$ along $j(M)$ in the above definitions was 
already well established in \cite{DG} and \cite{FG} for isometric 
bendings in both the local  and global situation. 
\medskip 

An infinitesimal bending $\tau$ of an isometric immersion 
$f\colon M^n\to\R^{n+p}$, $p\geq 2$, is called a \emph{genuine 
infinitesimal bending} if $\tau$ does not extend in the singular 
sense when restricted to any open subset of $M^n$. If $f$ 
admits such a bending we say that it is 
\emph{genuinely infinitesimally bendable}. 
As one expects, trivial infinitesimal bending are never genuine. 
If $f(M)\subset\R^{n+\ell}\subset\R^{n+p}$, $\ell<p$, and 
$e\in\R^{n+p}$ is orthogonal to $\R^{n+\ell}$ then $\tau=\phi e$ 
for $\phi\in C^\infty(M)$ is another example of an infinitesimal 
bending that is not genuine.
\medskip

Recall that an isometric immersion  $f\colon M^n\to\R^{n+p}$ is 
said to be \emph{$r$-ruled} if there exists an $r$-dimensional 
smooth totally geodesic tangent distribution whose leaves (rulings) 
are mapped diffeomorphically by $f$ to open subsets of affine 
subspaces of $\R^{n+p}$. 

\begin{theorem}\label{local1}
Let $f\colon M^n\to\R^{n+p}$, $n>2p\geq 4$, be an isometric immersion
and let $\tau$ be an infinitesimal bending of $f$. Then along each 
connected component of an open and dense subset either $\tau$ extends 
in the singular sense or $f$ is $r$-ruled with $r\geq n-2p$. 
\end{theorem}

The following is an immediate consequence of the above result.

\begin{corollary}
Let $f\colon M^n\to\R^{n+p}$, $n>2p\geq 4$, be a genuinely 
infinitesimally bendable isometric immersion. Then $f$ is 
$r$-ruled with $r\geq n-2p$ along connected components of an 
open dense subset of $M^n$. 
\end{corollary}

We say that $f\colon M^n\to\R^{n+p}$ is \emph{genuinely infinitesimally 
rigid} if given any infinitesimal bending $\tau$ of $f$ there is an open
dense subset of $M^n$ such that $\tau$ restricted to any connected 
component extends in the singular sense.
\vspace{1ex}

Theorem \ref{local1} also has the following two consequences.

\begin{corollary}
Let $f\colon M^n\to\R^{n+p}$, $n>2p\geq 4$, be an isometric 
immersion. If $M^n$ has positive Ricci curvature then $f$ is 
genuinely infinitesimally rigid. 
\end{corollary}

\begin{corollary}
Let $g\colon M^n\to\Sf^{n+p-1}$, $n>2p\geq 4$, be an isometric immersion
and let $f=i\circ g$ where $i\colon\Sf^{n+p-1}\to\R^{n+p}$ denotes the
umbilical inclusion. Then $f$ is genuinely infinitesimally rigid. 
\end{corollary}

A special class of ruled submanifolds are the ones with a relative 
nullity foliation. The \emph{relative nullity} subspace $\Delta(x)$ 
of $f\colon M^n\to\R^{n+p}$ at $x\in M^n$ is the kernel of the second 
fundamental form $\a\colon TM\times TM\to N_fM$ with values in the 
normal bundle, that is,
$$
\Delta(x)=\{X\in T_xM: \a(X,Y)=0\;\;\mbox{for all}\;\;Y\in T_xM\}.
$$
The dimension $\nu(x)$ of $\Delta(x)$ is called the \emph{index of 
relative nullity} of $f$ at $x\in M^n$. It is a standard fact that 
the submanifold is ruled by the leaves of the relative nullity 
distribution on any open subset of $M^n$ where the index of 
relative nullity $\nu>0$ is constant.
\medskip

In the case of low codimension, with a substantial additional effort 
we obtain a better lower bound for the dimension of the rulings. 

\begin{theorem}\label{local2}
Let $f\colon M^n\to\R^{n+p}$, $n>2p$, be a genuinely infinitesimally 
bendable isometric immersion. If $2\leq p\leq 5$, then one of the 
following holds along any connected component, say $U$, of an open 
dense subset of $M^n$:
\begin{itemize}
\item[(i)] $f|_U$ is $\nu$-ruled by leaves of relative 
nullity with $\nu\geq n-2p$.
\item [(ii)] $f|_U$ has $\nu<n-2p$ at any point and is $r$-ruled
with $r\geq n-2p+3$.
\end{itemize}
\end{theorem}

For $p=2$ notice that we are always in case $(i)$ since a 
$(n-1)$-ruled submanifold in that codimension has index of 
relative nullity $\nu\geq n-3$ at any point.\vspace{1ex} 
 
Dajczer and Gromoll \cite{DG} proved that along connected components 
of an open dense subset an isometrically deformable compact Euclidean 
submanifold of dimension at least five and codimension two is either
isometrically rigid or is contained in a deformable hypersurface 
(with possible singularities) and any isometric deformation of the 
former is given by an isometric deformation of the latter. This result 
was extended by Florit and Guimarães \cite{FG} to other low codimensions.  
The next result of similar nature concerns infinitesimal bendings of 
submanifolds in codimension two. 

\begin{theorem}\label{maincompact}
Let $f\colon M^n\to\R^{n+2}$, $n\geq 5$, be an isometric immersion   
of a compact Riemannian manifold with no open flat subset. For
any infinitesimal bending $\tau$ of $f$ one of the following 
holds along any connected component, say $U$, of an open dense 
subset of $M^n$:
\begin{itemize}
\item[(i)] The infinitesimal bending $\tau|_U$ extends in the 
singular sense.
\item[(ii)] There is an orthogonal splitting 
$\R^{n+2}=\R^{n+1}\oplus\spa\{e\}$ so that
$f(U)\subset\R^{n+1}$ and $\tau|_U=\tau_1+\tau_2$ is a sum of 
infinitesimal bendings that extend in the singular sense where 
$\tau_1\in\R^{n+1}$ and $\tau_2=\phi e$ for $\phi\in C^\infty(U)$. 
\end{itemize}
\end{theorem}

It follows from the proof  that the assumption on the open flat subset 
can be replaced by the weaker hypothesis that there is no open subset 
of $M^n$ where the index of relative nullity satisfies $\nu\geq n-1$. 
Moreover, we will see that cases $(i)$ and $(ii)$ are not disjoint.
\vspace{1ex}

In the last section of the paper, we discuss why the local results 
given above also hold if the ambient space is a nonflat space form.

\section{The associated tensor}

In this section, we discuss several properties of a tensor
associated to an infinitesimal bending  called in the classical 
theory of surfaces the associated rotation field; for instance
see \cite{Sp}. For basic facts on infinitesimal bendings we 
refer to \cite{DR}, \cite{DT}, \cite{DV} and \cite{GR}.
\medskip

In the sequel, let $\tau$ denote an infinitesimal bending of a 
isometric immersion $f\colon M^n\to\R^{n+p}$.  Then the section 
$L\in\Gamma(\Hom(TM,f^*T\R^{n+p}))$ is the tensor defined as 
$$
LX=\nab_X\tau
$$
where $\nab$ is the Levi-Civita connection in $\R^{n+p}$.
Hence \eqref{infbend} can be written as 
\be\label{skew}
\<LX,f_*Y\>+\<LY,f_*X\>=0
\ee
for any $X,Y\in\mathfrak{X}(M)$. 

Let $B\colon TM\times TM\to f^*T\R^{n+p}$ the symmetric tensor   
defined as
$$
B(X,Y)=(\nab_XL)Y
$$
for any $X,Y\in\mathfrak{X}(M)$. If $\tau$ is an immersion notice 
that $B$ is nothing else than its second fundamental form. 

\begin{proposition} The tensor $B$ satisfies
\be\label{segderL}
(\nab_XB)(Y,Z)-(\nab_YB)(X,Z)=-LR(X,Y)Z
\ee
for any $X,Y,Z\in\mathfrak{X}(M)$. 
\end{proposition}
\proof Use that
\be\label{form}
(\nab_XB)(Y,Z)
=\nab_X(\nab_YL)Z-(\nab_{\nabla_XY}L)Z-(\nab_Y L)\nabla_XZ
\ee
and the definition of the curvature tensor.\vspace{2ex}\qed

The metrics $g_t$ induced by $f_t=f+t\tau$ satisfy
\be\label{der metrica}
\d/\d t|_{t=0}g_t(X,Y)=0
\ee
for any $X,Y\in\mathfrak{X}(M)$. Hence, the Levi-Civita 
connections and curvature tensors of $g_t$ verify
\be\label{der conex}
\d/\d t|_{t=0}\nabla^{t}_X Y=0 
\ee
and
\be\label{der curv}
\d/\d t|_{t=0}g_t(R^t(X,Y)Z,W)=0
\ee 
for any $X,Y,Z,W\in\mathfrak{X}(M)$.
Taking the derivative with respect to $t$ at $t=0$ of the
Gauss formula for  $f_t$, namely, of
$$
\nab_Xf_{t*}Y=f_{t*}\nabla^t_XY+\a^t(X,Y),
$$
we obtain
\be\label{derL}
B(X,Y)=\d/\d t|_{t=0}\alpha^t(X,Y).
\ee

Taking tangent and normal components with respect to $f$ we have
$$
B(X,Y)=f_*\Y(X,Y)+\beta(X,Y)
$$
where the  tensors $\Y\colon TM\times TM\to TM$ and 
$\beta\colon TM\times TM\to N_fM$ are also symmetric. 

\begin{proposition} The tensor $\Y\colon TM\times TM\to TM$ 
satisfies
\be\label{parte}
\<\alpha(X,Y),LZ\>+\<\Y(X,Y),Z\>=0
\ee
for any $X,Y,Z\in\mathfrak{X}(M)$.
\end{proposition}

\proof Given $\eta(t)\in\Gamma(N_{f_t}M)$, let 
$\Y_\eta$ be the tangent vector field given by
$$
f_*\Y_\eta=(\d/\d t|_{t=0}\eta(t))_{f_*TM}.
$$
The derivative of $\<f_{t*}Z,\eta(t)\>=0$ with respect to $t$ at 
$t=0$ yields
$$
\<\eta,LZ\>+\<\Y_{\eta},Z\>=0
$$
where $Z\in\mathfrak{X}(M)$ and $\eta=\eta(0)$. In particular,
$$
\<\alpha(X,Y),LZ\>+\<\Y_{\alpha(X,Y)},Z\>=0
$$
for any $X,Y,Z\in\mathfrak{X}(M)$. On the other hand, we obtain from 
\eqref{derL} that
$$
\Y_{\alpha(X,Y)}=\Y(X,Y)
$$
for any $X,Y\in\mathfrak{X}(M)$.\qed

\begin{proposition} The tensor $\beta\colon TM\times TM\to N_fM$ 
satisfies
\begin{align}
\label{derGauss}
\<\beta(X,W),\alpha(Y,Z)\>\,+\,&\<\alpha(X,W),\beta(Y,Z)\>\nonumber\\
&=\<\beta(X,Z),\alpha(Y,W)\>+\<\alpha(X,Z),\beta(Y,W)\>
\end{align}
and
\begin{align}\label{casicodazzi}
(\nabla^{\perp}_X\beta)(Y,Z)-&\,(\nabla_Y^{\perp}\beta)(X,Z)\nonumber\\
&=\,\alpha(Y,\Y(X,Z))-\alpha(X,\Y(Y,Z))-(LR(X,Y)Z)_{N_fM}
\end{align}
for any $X,Y,Z,W\in\mathfrak{X}(M)$.
\end{proposition}

\proof To prove \eqref{derGauss} take the derivative with respect to 
$t$ at $t=0$ of the Gauss equations for $f_t$, that is, of
$$
g_t(R^t(X,Y)Z,W)=g_t(\a^t(X,W),\a^t(Y,Z))-g_t(\a^t(X,Z),\a^t(Y,W))
$$
and use \eqref{der metrica}, \eqref{der curv} and \eqref{derL}.

Using \eqref{form} we have
$$
((\nab_XB)(Y,Z))_{N_fM}=
\alpha(X,\Y(Y,Z))+(\nabla_X^{\perp}\beta)(Y,Z)
$$
and \eqref{casicodazzi} follows from 
\eqref{segderL}.\vspace{2ex}\qed

We discuss next the simplest examples of infinitesimal bendings.

\begin{examples} \label{examples}\emph{
$(1)$ If $\tau$ is a trivial infinitesimal bending of 
$f\colon M^n\to\R^{n+p}$, $p\geq 2$, then we have from
the references that
$$
\tau=\mathcal{D}f(x)+w
$$
where $\mathcal{D}$ is a skew-symmetric linear transformation of 
$\R^{n+p}$ and $w\in\R^{n+p}$. Take $\lambda\in\Gamma(f^*T\R^{n+p})$ 
such that 
$F\colon \tilde{M}^{n+1}=M^n\times(-\epsilon,\epsilon)\to \R^{n+p}$, 
given by $F(x,t)=f(x)+t\lambda(x)$, is an immersion for $t\neq 0$. 
Then $\tau$ extends in the singular sense since
$$
\tilde{\tau}(x,t)=\tau+t\mathcal{D}\lambda
$$ 
is a (trivial) infinitesimal bending of $F$ on the open subset where $F$ is an immersion.\vspace{1ex}\\
$(2)$ The first normal space of 
$f\colon M^n\to\R^{n+p}$ at $x\in M^n$ is 
$$
N_1(x)=\spa\{\a(X,Y):X,Y\in T_xM\}.
$$
Then $\tau=f_*Z+\delta$ is an infinitesimal bending if 
$Z\in\mathfrak{X}(M)$ is a Killing field and 
$\delta\in \Gamma(N_1^\perp)$ is a smooth normal vector field.
}\end{examples}

\section{Flat bilinear forms}

Flat bilinear forms were introduced by J. D.  Moore \cite{Mo} after 
the pioneering work of E. Cartan to deal with rigidity questions on 
isometric immersions in space forms. In this paper, it is shown 
that they are also very helpful in the study of similar questions 
for infinitesimal bendings of submanifolds.
\vspace{1ex}

Let $V$ and $U$ be finite dimensional real vector spaces and 
let $W^{p,q}$ be a real vector space of dimension $p+q$
endowed with an indefinite inner product of type $(p,q)$. 
A bilinear form $\mathcal{B}\colon V\times U\to W^{p,q}$ is said 
to be \emph{flat} if 
$$
\<\mathcal{B}(X,Z),\mathcal{B}(Y,W)\>
-\<\mathcal{B}(X,W),\mathcal{B}(Y,Z)\>=0
$$
for all $X,Y\in V$ and $W,Z\in U$.  
Then $X\in V$ is called a (left) \emph{regular element} of 
$\mathcal{B}$ if  
$$
\dim\mathcal{B}_X(U)=\max\{\dim\mathcal{B}_Y(U)\colon Y\in V\}
$$
where $\mathcal{B}_X(Y)=\mathcal{B}(X,Y)$ for any $Y\in U$. 
The set $RE(\mathcal{B})$ of regular elements of $\mathcal{B}$
is open dense in $V$.
\medskip

The following basic fact was given in \cite{Mo}. 

\begin{lemma}\label{fbn}
Let $\mathcal{B}\colon V\times U\to W$ be a flat bilinear form. 
If $Y\in RE (\mathcal{B})$ then 
$$
\mathcal{\mathcal{B}}(X,\ker\mathcal{B}_Y)
\subset \mathcal{B}_Y(U)\cap\mathcal{B}_Y(U)^{\perp}
$$
for any $X\in V$.
\end{lemma}

The next is a fundamental result in the theory of symmetric flat
bilinear forms. It turns out to be 
false for $p\geq 6$ as shown in \cite{DF2}.

\begin{lemma}\label{main} 
Let $\mathcal{B}\colon V^n\times V^n\to W^{p,q}$, $p\leq 5$ 
and $p+q<n$, be a symmetric flat bilinear form and set 
$$
\mathcal{N}(\mathcal{B})=\{X\in V:
\mathcal{B}(X,Y)=0\;\;\mbox{for all}\;\;Y\in V\}.
$$
If $\dim \mathcal{N}(\mathcal{B})\leq n-p-q-1$ then there is 
an orthogonal decomposition
$$
W^{p,q}=W_1^{\ell,\ell}\oplus W_2^{p-\ell,q-\ell},\; 1\leq\ell\leq p,
$$
such that the $W_j$-components $\mathcal{B}_j$ of $\mathcal{B}$ 
satisfy:
\begin{itemize}
\item[(i)] $\mathcal{B}_1$ is nonzero and
$$
\<\mathcal{B}_1(X,Y),\mathcal{B}_1(Z,W)\>=0
$$
for all $X,Y,Z,W\in V$.
\item[(ii)] $\mathcal{B}_2$ is flat and 
$\dim\mathcal{N}(\mathcal{B}_2)\geq n-p-q+2\ell$.
\end{itemize}
\end{lemma}

\proof See \cite{DF} or \cite{DT}. \qed

\section{The local results}

In this section we give the proofs the local theorems 
in the introduction. A key ingredient is the following
result due to Florit and Guimarães \cite{FG}.

\begin{proposition}\label{nowhere} Let $f\colon M^n\to\R^{n+p}$ 
be an isometric immersion and let $D$ be a smooth tangent 
distribution of dimension $d>0$. Assume that there does not exist 
an open subset $U\subset M^n$ and $Z\in\Gamma(D|_U)$ such 
that the map $F\colon U\times\R\to\R^{n+p}$ given by
$$
F(x,t)=f(x)+tf_*Z(x)
$$
is a singular extension of $f$ on some open neighborhood of 
$U\times\{0\}$. Then for any $x\in M^n$ there is an open 
neighborhood $V$ of the origin in $D(x)$ such that 
$f_*(x)V\subset f(M)$. Hence $f$ is $d$-ruled along each 
connected component of an open dense subset of $M^n$.
\end{proposition}

\proof See \cite{DT} or \cite{FG}.\qed

\subsection{The first local result}

We first associate to an infinitesimal bending a flat bilinear form.

\begin{lemma}\label{thetaflat} 
Let $\tau$ be an infinitesimal bending of an isometric 
immersion $f\colon M^n\to\R^{n+p}$. Then the bilinear form
$\theta\colon TM\times TM\to N_fM\oplus N_fM$ defined at 
any point of $M^n$ by
\be\label{theta}
\theta(X,Y)=(\alpha(X,Y)+\beta(X,Y),\alpha(X,Y)-\beta(X,Y)) 
\ee
is flat  with respect to the inner product in 
$N_fM\oplus N_fM$  given by
$$
\<\!\<(\xi_1,\eta_1),(\xi_2,\eta_2)\>\!\>_{N_fM\oplus N_fM}
=\<\xi_1,\xi_2\>_{N_fM}-\<\eta_1,\eta_2\>_{N_fM}.
$$
\end{lemma}

\proof A straightforward computation shows that
\begin{align*}
&\frac{1}{2}\left(\<\!\<\theta(X,Z),\theta(Y,W)\>\!\>
-\<\!\<\theta(X,W),\theta(Y,Z)\>\!\>\right)
=\<\beta(X,Z),\alpha(Y,W)\>\\
&\;\;\;+\<\alpha(X,Z),\beta(Y,W)\> 
-\<\beta(X,W),\alpha(Y,Z)\>-\<\alpha(X,W),\beta(Y,Z)\>,
\end{align*}
and the proof follows from \eqref{derGauss}.\vspace{2ex}\qed

An isometric immersion $f\colon M^n\to\R^{n+p}$ is called
\emph{$1$-regular} if the first normal spaces $N_1(x)$ have 
constant dimension $k\leq p$ on $M^n$ and thus form a subbundle
$N_1$ of rank $k$ of the normal bundle. Under the $1$-regularity 
assumption we have the following equivalent statement.

\begin{lemma}\label{cor}  Assume that $f$ 
is $1$-regular and let $\beta_1\colon TM\times TM\to N_1$ be
the $N_1$-component of $\beta$. Then the bilinear form
$\hat\theta\colon TM\times TM\to N_1\oplus N_1$ defined at any 
point by
\be\label{deftheta1}
\hat\theta(X,Y)
=\left(\alpha(X,Y)+\beta_1(X,Y),\alpha(X,Y)-\beta_1(X,Y)\right) 
\ee
is flat with respect to the inner product induced on 
$N_1\oplus N_1$. 
\end{lemma}

\noindent\emph{Proof of Theorem \ref{local1}:}
Let $\tau$ be an infinitesimal bending of $f$. With the 
use of \eqref{skew} and \eqref{parte} we easily obtain
\be\label{above}
\<f_*X+\nab_XY,LX+\nab_X LY\>=\<\alpha(X,Y),\beta(X,Y)\>
\ee
for any $X,Y\in\mathfrak{X}(M)$.

By Lemma \ref{thetaflat} we have at any point of $M^n$  that 
the symmetric tensor $\theta$ is flat. Given $Y\in RE(\theta)$ at 
a point denote $D=\ker\theta_Y$ where $\theta_Y(X)=\theta(Y,X)$. 
Notice that $Z\in D$ means that $\a(Y,Z)=0=\beta(Y,Z)$.

Let $U\subset M^n$ be an open subset where $Y\in\mathfrak{X}(U)$
satisfies $Y\in RE(\theta)$ and $D$ has dimension $d$ at any 
point. Lemma \ref{fbn} gives
$$
\<\!\<\theta(X,Z),\theta(X,Z)\>\!\>=0
$$
for any $X\in\mathfrak{X}(U)$ and $Z\in\Gamma(D)$. Equivalently, 
the right hand side of \eqref{above} vanishes and thus
\be\label{vanish}
\<f_*X+\nab_XZ,LX+\nab_X LZ\>=0
\ee
for any $X\in\mathfrak{X}(U)$ and $Z\in\Gamma(D)$. 

Assume that there exists a nowhere vanishing $Z\in\Gamma(D)$ 
defined on an open subset $V$ of $U$ such that
$F\colon V\times(-\epsilon,\epsilon)\to\R^{n+p}$ given by
$$
F(x,t)=f(x)+tf_*Z(x)
$$
is a singular extension of $f|_{V}$. The map
$\tilde{\tau}\colon V\times(-\epsilon,\epsilon)\to\R^{n+p}$ 
given by
$$
\tilde{\tau}(x,t)=\tau(x)+tLZ(x)
$$
is an infinitesimal bending as well as an extension of 
$\tau|_{V}$ in the singular sense. In fact, 
$$
\<F_*\d_t,\nab_{\d_t}\tilde{\tau}\>=\<f_*Z,LZ\>=0,
$$
$$
\<\nab_{\d_t}\tilde{\tau},F_*X\>+\<\nab_X\tilde{\tau},F_*\d_t\>
=\<LZ,f_*X+t\nab_XZ\>+\<LX+t\nab_XLZ,f_*Z\>=0
$$
and
$$
\<F_*X,\nab_X\tilde{\tau}\>=\<f_*X+t\nab_XZ,LX+t\nab_XLZ\>=0
$$
where the last equality follows from \eqref{vanish}.

Let $W\subset U$ be an open subset such that a $Z\in\Gamma(D)$ as 
above does not exist along any open subset of $W$. By 
Proposition \ref{nowhere} the immersion is $d$-ruled
along any connected component of an open dense subset of $W$. 
Moreover, we have
$d=\dim D=n-\dim\Ima(\theta_Y)\geq n-2p$. \qed

\begin{remark} {\em In Theorem \ref{local1} if $f$ 
is $1$-regular with $\dim N_1=q<p$ we obtain the better lower
bound $r\geq n-2q$ since the proof still works making use of 
Lemma \ref{cor} instead of Lemma \ref{thetaflat}.
}\end{remark}

\subsection{The second local result}

Let $F\colon\tilde{M}^{n+1}\to\R^{n+p}$ be an isometric immersion
and let $\tilde{\tau}$ be an infinitesimal bending of $F$.
Given an isometric embedding $j\colon M^n\to\tilde{M}^{n+1}$  
consider the composition of isometric immersions 
$f=F\circ j\colon M^n\to\R^{n+p}$. Then
$\tau=\tilde{\tau}|_{j(M)}$ is an infinitesimal bending of $f$. 
It is easy to see that 
$$
B(X,Y)=\tilde{B}(X,Y)+\<\nab_XY,F_*\eta\>\tilde{L}\eta
$$ 
for $\eta\in\Gamma(N_jM)$ of unit length and $X,Y\in\mathfrak{X}(M)$. 
Then \eqref{parte} gives
$$
\<\beta(X,Y),F_*\eta\>+\<\alpha^f(X,Y),\tilde{L}\eta\>=0
$$
for any $X,Y\in\mathfrak{X}(M)$. We will see that satisfying a 
condition of this type may guarantee that an infinitesimal bending 
is not genuine. In fact, this was already proved by Florit \cite{Fl} 
in a special case.
\medskip

We say that an infinitesimal bending of a given isometric immersion 
$f\colon M^n\to\R^{n+p}$, $p\geq 2$, satisfies the \emph{condition $(*)$} 
if there is $\eta\in\Gamma(N_fM)$ nowhere vanishing and $\xi\in\Gamma(R)$, 
where $R$ is determined by the orthogonal splitting $N_fM=P\oplus R$ 
and $P=\spa\{\eta\}$, such that
\be\label{criterio}
B_\eta+A_\xi=0
\ee
where $B_\eta=\<\beta,\eta\>$. We choose $\eta$
of unit length for simplicity.
Thus, that \eqref{criterio} holds means 
\be\label{criterios}
\<\beta(X,Y),\eta\>+\<\alpha(X,Y),\xi\>=0
\ee
for any $X,Y\in\mathfrak{X}(M)$.
\medskip

The following result is of independent interest since it does not 
require the codimension to satisfy $p\leq 5$ as is the case in 
Theorem \ref{local2}.

\begin{theorem}\label{local}
Let $f\colon M^n\to\R^{n+p}$, $p\geq 2$, be an isometric immersion 
and let $\tau$ be an infinitesimal bending of $f$ that satisfies 
the condition $(*)$. Then along each connected component of an 
open and dense subset of $M^n$ either $\tau$ extends in the singular 
sense or $f$ is $r$-ruled with $r\geq n-2p+3$.
\end{theorem}

As before there is the following immediate consequence.

\begin{corollary}\label{localgen}
Let $f\colon M^n\to\R^{n+p}$, $p\geq 2$, be an isometric immersion 
and let $\tau$ be a genuine infinitesimal bending of $f$ that 
satisfies the condition $(*)$. Then $f$ is $r$-ruled with 
$r\geq n-2p+3$ on connected components of an open dense subset 
of $M^n$.
\end{corollary}

When $\tau$ satisfies the condition $(*)$ we may extend the tensor 
$L$ to the tensor $\bar{L}\in\Gamma(\Hom(TM\oplus P,f^*T\R^{n+p})$ 
by defining
$$
\bar{L}\eta=f_*Y+\xi
$$
where $Y\in\mathfrak{X}(M)$ is given by
$$
\<Y,X\>+\<LX,\eta\>=0
$$
for any $X\in\mathfrak{X}(M)$. Then $\bar{L}$ satisfies 
$$
\<\bar{L}X,\eta\>+\<f_*X,\bar{L}\eta\>=0
$$
for any $X\in\mathfrak{X}(M)$.

Given $\lambda\in\Gamma (f_*TU\oplus P)$ nowhere vanishing where 
$U$ is an open subset of $M^n$, we define the map
$F\colon U\times (-\epsilon,\epsilon)\to\R^{n+p}$
by
\be\label{F}
F(x,t)=f(x)+t\lambda(x).
\ee
Notice that $F$ is not an immersion at least for $t=0$ at points 
where $\lambda$ is tangent to $U$. Then let 
$\tilde{\tau}\colon U\times (-\epsilon,\epsilon)\to\R^{n+p}$
be the map given by
\be\label{tilde}
\tilde{\tau}(x,t)=\tau(x)+t\bar{L}\lambda(x). 
\ee
We have
$$
\<F_*\d_t,\nab_{\d_t}\tilde{\tau}\>=0.
$$
Moreover, since $\<\bar{L}\lambda,\lambda\>=0$ we obtain
$$
\<\nab_{\d_t}\tilde{\tau},F_*X\>+\<\nab_X\tilde{\tau},F_*\d_t\>
=\<\bar{L}\lambda,f_*X\>
+\<LX,\lambda\>+tX\<\bar{L}\lambda,\lambda\>=0
$$
for any $X\in\mathfrak{X}(M)$ and $t\in (-\epsilon,\epsilon)$.
Thus $\tilde{\tau}$ is an infinitesimal bending 
of $F$ on the open subset $\tilde{U}$ of $U\times (-\epsilon,\epsilon)$
where $F$ is an immersion if and only if 
$$
\<F_*X,\nab_X \tilde{\tau}\>=0,
$$ 
or equivalently, if
$$
\<f_*X+t\nab_X\lambda,LX+t\nab_X\bar{L}\lambda\>=0
$$
for any $X\in\mathfrak{X}(M)$.

In the sequel we take  $F$ restricted to $\tilde{U}$.
By the above, in order to have that $\tilde{\tau}$ is an infinitesimal 
bending of $F$ the strategy is to make use of the condition $(*)$ to 
construct a subbundle $D\subset f_*TM\oplus P$ such that 
$$
\<f_*X+\nab_X \lambda,LX+\nab_X\bar{L}\lambda\>=0
$$
for any $X\in\mathfrak{X}(M)$ and any $\lambda\in\Gamma(D)$.

\begin{lemma} Assume that $\tau$ satisfies the condition $(*)$.
Then
\be\label{impext}
\<f_*X+\nab_X \lambda,LX+\nab_X\bar{L}\lambda\>
=\<(\nab_X\lambda)_R,(\nab_X \bar{L})\lambda\>
\ee
where $X\in\mathfrak{X}(M)$, $\lambda\in\Gamma(f_*TM\oplus P)$ and
$$
(\nab_X\bar{L})\lambda=\nab_X\bar{L}\lambda-\bar{L}\nabla'_X\lambda,
$$ 
being $\nabla'$ the connection induced on $f_*TM\oplus P$.
\end{lemma}

\proof Set $\lambda=f_*Z+\phi\eta$ where $Z\in\mathfrak{X}(M)$
and $\phi\in C^\infty(M)$. Then 
\begin{align}\label{eq 0}
&\<f_*X+\nab_X \lambda,LX+\nab_X\bar{L}\lambda\>=
\<f_*(\nab_X\lambda)_{TM}+(\nab_X\lambda)_P
+(\nab_X\lambda)_R,\nab_X\bar{L}\lambda\>\nonumber\\
&\;\;+\<\nab_X\lambda,LX\>+\<f_*X,\nab_X\bar{L}\lambda\>\nonumber\\
&=\<f_*(\nab_X\lambda)_{TM},(\nab_XL)Z+L\nabla_X Z
+X(\phi)\bar{L}\eta+\phi\nab_X\bar{L}\eta\>\nonumber\\
&\;\;+(\<A_\eta X,Z\>+X(\phi))\<\eta,(\nab_XL)Z+L\nabla_XZ
+X(\phi)\bar{L}\eta+\phi\nab_X\bar{L}\eta\>\nonumber\\
&\;\;+\<(\nab_X\lambda)_R,\nab_X\bar{L}\lambda\>
+\<\nab_X\lambda,LX\>+\<f_*X,\nab_X\bar{L}\lambda\>
\end{align}
for any $X\in\mathfrak{X}(M)$. Using \eqref{parte} and 
\eqref{criterios} we obtain 
\begin{align}\label{eq 1}
\<(\nab_X\lambda)_{TM},&(\nab_XL)Z+L\nabla_XZ\>\nonumber\\
&=-\<L(\nab_X\lambda)_{TM},\alpha(X,Z)\>
-\phi\<A_\eta X,L\nabla_X Z\>
\end{align}
and
\begin{align}\label{eq 2}
\<(\nab_X\lambda)_{TM},&X(\phi)\bar{L}\eta
+\phi\nab_X\bar{L}\eta\>=\phi\<(\nab_X\lambda)_{TM},\nabla_XY\>\nonumber\\
\!&
-X(\phi)\<L(\nab_X\lambda)_{TM},\eta\>
-\phi\<\alpha(X,(\nab_X\lambda)_{TM}),\xi\> 
\end{align}
where for the first term in the right hand side of \eqref{eq 2} we have 
\begin{align}\label{eq 5}
\<(\nab_X\lambda)_{TM},\nabla_XY\>
=&\,X\<(\nab_X\lambda)_{TM},Y\>
-\<\nabla_X(\nab_X\lambda)_{TM},Y\>\nonumber\\
=&-X\<L(\nab_X\lambda)_{TM}),\eta\>
+\<L\nabla_X(\nab_X\lambda)_{TM},\eta\>\nonumber\\ 
=&-\<(\nab_XL)(\nab_X\lambda)_{TM},\eta\>
-\<L(\nab_X\lambda)_{TM},\nab_X\eta\>\nonumber\\
=&\,\<\alpha(X,(\nab_X\lambda)_{TM}),\xi\>
-\<L(\nab_X\lambda)_{TM},\nab_X\eta\>.
\end{align}
Moreover,
\begin{align}\label{eq 3}
\<\eta,(\nab_XL)Z+L\nabla_XZ\>=&
-\<\alpha(X,Z),\xi\>+\<\eta,L\nabla_XZ\>,
\end{align}
\begin{align}\label{eq 4}
\<\eta,X(\phi)\bar{L}\eta+\phi\nab_X\bar{L}\eta\>=&
-\phi\<\nab_X\eta,\bar{L}\eta\>\nonumber\\
=&-\phi\<LA_\eta X,\eta\>-\phi\<\nabla^\perp_X\eta,\xi\>
\end{align}
and 
\begin{align}\label{eq 6}
&\<\nab_X\lambda,LX\>+\<f_*X,\nab_X\bar{L}\lambda\>
=-\<\nab_X X,\bar{L}\lambda\>-\<\lambda,\nab_XLX\>\nonumber\\
=&-\<\nabla_XX,\bar{L}\lambda\>-\<\alpha(X,X),\bar{L}\lambda\>
-\<\lambda,L\nabla_X X\>-\<\lambda,(\nab_XL)X\>=0.
\end{align}
Now a straightforward computation replacing \eqref{eq 1} through 
\eqref{eq 6} in \eqref{eq 0} yields
\begin{align*}
&\<f_*X+\nab_X \lambda,LX+\nab_X\bar{L}\lambda\>=
\<(\nab_X\lambda)_R,\nab_X\bar{L}\lambda\>
-\<L(\nab_X\lambda)_{TM},\alpha(X,Z)_R\>\\
&-\phi\<L(\nab_X\lambda)_{TM},\nabla^\perp_X\eta\>
-\<\alpha(X,Z),\bar{L}(\nab_X\lambda)_P\>
-\phi\<\nabla^\perp_X\eta,\bar{L}(\nab_X\lambda)_P\>\\
&=\,\<(\nab_X\lambda)_R,(\nab_X\bar{L})\lambda\>.\qed
\end{align*}

In view of \eqref{impext} the next step is to construct a subbundle 
$D\subset f_*TM\oplus P$ such that 
\be\label{requisito}
\<(\nab_X\lambda)_R,(\nab_X\bar{L})\lambda\>=0
\ee
for any $X\in\mathfrak{X}(M)$ and $\lambda\in\Gamma(D)$. 

\begin{lemma}\label{varphi} Assume that $\tau$ satisfies the 
condition $(*)$. Then the bilinear form 
$\varphi\colon TM\times f_*TM\oplus P\to R\oplus R$ defined~by
$$
\varphi(X,\lambda)
=((\nab_X\lambda)_R+((\nab_X\bar{L})\lambda)_R,
(\nab_X\lambda)_R-((\nab_X\bar{L})\lambda)_R)
$$
is flat with respect to the indefinite inner product given by
$$
\<\!\<(\xi_1,\mu_1),(\xi_2,\mu_2)\>\!\>_{R\oplus R}
=\<\xi_1,\xi_2\>_R-\<\mu_1,\mu_2\>_R.
$$
\end{lemma}

\proof We need to show that
$$
\Theta=\<\!\<\varphi(X,\lambda),\varphi(Y,\delta)\>\!\>
-\<\!\<\varphi(X,\delta),\varphi(Y,\lambda)\>\!\>=0
$$ 
for any $X,Y\in\mathfrak{X}(M)$ and $\lambda,\delta\in f_*TM\oplus P$.
We have 
\begin{align*}
\frac{1}{2}\Theta=&\<(\nab_X\lambda)_R,((\nab_Y\bar{L})\delta)_R\>
+\<(\nab_Y\delta)_R,((\nab_X\bar{L})\lambda)_R\>\nonumber\\
&\;\;-\<(\nab_X\delta)_R,((\nab_Y\bar{L})\lambda)_R\>
-\<(\nab_Y\lambda)_R,((\nab_X\bar{L})\delta)_R\>.
\end{align*}
Clearly $\Theta=0$ if $\lambda,\delta\in\Gamma(P)$. 
If $\lambda,\delta\in\mathfrak{X}(M)$, then
\begin{align*}
\frac{1}{2}\Theta
=&\,\<\alpha(X,\lambda)_R,((\nab_Y\bar{L})\delta)_R\>
+\<\alpha(Y,\delta)_R,((\nab_X\bar{L})\lambda)_R\>\\
&-\<\alpha(X,\delta)_R,((\nab_Y\bar{L})\lambda)_R\>
-\<\alpha(Y,\lambda)_R,((\nab_X\bar{L})\delta)_R\>\\
=&\,\<\alpha(X,\lambda)_R,((\nab_Y L)\delta)_R\>
-\<A_\eta Y,\delta\>\<\alpha(X,\lambda)_R,\bar{L}\eta\>\\
&+\<\alpha(Y,\delta)_R,((\nab_X L)\lambda)_R\>
-\<A_\eta X,\lambda\>\<\alpha(Y,\delta)_R,\bar{L}\eta\>\\
&-\<\alpha(X,\delta)_R,((\nab_Y L)\lambda)_R\>
+\<A_\eta Y,\lambda\>\<\alpha(X,\delta)_R,\bar{L}\eta\>\\
&-\<\alpha(Y,\lambda)_R,((\nab_X L)\delta)_R\>
+\<A_\eta X,\delta\>\<\alpha(Y,\lambda)_R,\bar{L}\eta\>. 
\end{align*}
Using first \eqref{criterios} and then \eqref{derGauss}
we obtain
\begin{align*}
\frac{1}{2}\Theta=&\,\<\alpha(X,\lambda),\beta(Y,\delta)\>
+\<\alpha(Y,\delta),\beta(X,\lambda)\>\\
&-\<\alpha(X,\delta),\beta(Y,\lambda)\>
-\<\alpha(Y,\lambda),\beta(X,\delta)\>=0.
\end{align*}
Finally, we consider the case $\lambda=\eta$ and 
$\delta=Z\in\mathfrak{X}(M)$. Then
\begin{align*}
\frac{1}{2}\Theta=&\,\<\nabla^\perp_X\eta,((\nab_Y L)Z)_R\>
-\<A_\eta Y,Z\>\<\nabla^\perp_X\eta,\bar{L}\eta\>
+\<\alpha(Y,Z)_R,((\nab_X\bar{L})\eta)_R\>\\
-&\,\<\nabla^\perp_Y\eta,((\nab_X L)Z)_R\>
+\<A_\eta X,Z\>\<\nabla^\perp_Y\eta,\bar{L}\eta\>
-\<\alpha(X,Z)_R,((\nab_Y\bar{L})\eta)_R\>.
\end{align*}
Since
\begin{align*}
\<\nabla^\perp_X\eta,\bar{L}\eta\>
=&\,\<\nab_X\eta,\bar{L}\eta\>
+\<A_\eta X,\bar{L}\eta\>
=-\<\eta,\nab_X\bar{L}\eta\>-\<LA_\eta X,\eta\>\\
=&-\<\eta,(\nab_X\bar{L})\eta\>
\end{align*}
we obtain
\begin{align*}
\frac{1}{2}\Theta=&\,\<\nabla^\perp_X\eta,(\nab_Y L)Z\>
-\<\nabla^\perp_Y\eta,(\nab_X L)Z\>\\
&+\<\alpha(Y,Z),(\nab_X\bar{L})\eta\>
-\<\alpha(X,Z),(\nab_Y\bar{L})\eta\>.
\end{align*}
For the first term using \eqref{form}, \eqref{parte} and 
\eqref{criterios} we obtain
\begin{align*}
\<\nabla^\perp_X\eta, (\nab_Y L)Z\>=&\,X\<\eta,(\nab_YL)Z\>
-\<\eta,\nab_X(\nab_YL)Z\> +\<A_\eta X,(\nab_YL)Z\>\\
=&-X\<\alpha(Y,Z),\bar{L}\eta\>-\<\alpha(Y,Z),LA_\eta X\>\\
&-\<\eta,(\nab_X B)(Y,Z)+(\nab_{\nabla_X Y}L)Z
+(\nab_YL)\nabla_X Z\>\\ 
=&-\<(\nabla^\perp_X\alpha)(Y,Z)+\alpha(\nabla_X Y,Z)
+\alpha(Y,\nabla_X Z),\bar{L}\eta\>\\ 
&-\<\eta,(\nab_X B)(Y,Z)+(\nab_{\nabla_X Y}L)Z
+(\nab_YL)\nabla_X Z\>\\
&-\<\alpha(Y,Z),\nab_X\bar{L}\eta\>-\<\alpha(Y,Z),LA_\eta X\>
+\<A_{\alpha(Y,Z)}X,\bar{L}\eta\>\\
=&-\<(\nabla^\perp_X\alpha)(Y,Z),\bar{L}\eta\>
-\<\eta,(\nab_X B)(Y,Z)\>\\
&-\<\alpha(Y,Z),(\nab_X\bar{L})\eta\>-\<LA_{\alpha(Y,Z)}X,\eta\>.
\end{align*}
Likewise, we have
\begin{align*}
\<\nabla^\perp_Y\eta, (\nab_X L)Z\>=&
-\<(\nabla^\perp_Y\alpha)(X,Z),\bar{L}\eta\>
-\<\eta,(\nab_Y B)(X,Z)\>\\
&-\<\alpha(X,Z),(\nab_Y\bar{L})\eta\>-\<LA_{\alpha(X,Z)}Y,\eta\>.
\end{align*}
 From \eqref{segderL} and the Codazzi equation 
$$
(\nap_X\a)(Y,Z)=(\nap_Y\a)(X,Z)
$$
we obtain
\begin{align*}
\frac{1}{2}\Theta=\<L(R(X,Y)Z-A_{\alpha(Y,Z)}X+A_{\alpha(X,Z)}Y),\eta\>.
\end{align*}
Hence $\Theta=0$ from the Gauss equation.
\vspace{2ex}\qed

\noindent \emph{Proof of Theorem \ref{local}:} 
By Lemma \ref{varphi} there is a flat bilinear form $\varphi$. 
Let $U$ be an open subset of $M^n$ where there is 
$Y\in\mathfrak{X}(U)$ such that $Y\in RE(\varphi)$ and 
$D=\ker\varphi_Y$ has dimension $d$ at any point. 
Then Lemma \ref{fbn} gives 
$$
\<\!\<\varphi(X,\lambda),\varphi(X,\lambda)\>\!\>=0
$$
for any $X\in\mathfrak{X}(U)$ and $\lambda\in\Gamma(D)$. 
Notice that this implies that \eqref{requisito} holds for 
any $\lambda\in\Gamma(D)$. Whenever there is a nonvanishing
$\lambda\in\Gamma(D)$ on an open subset $V\subset U$ 
such that \eqref{F} defines a singular extension of $f|_V$, then 
$\tau|_V$ extends in the singular sense by means of \eqref{tilde}. 

Let $W\subset U$ be an open subset where  $\lambda\in\Gamma(D)$
as above does not exist along any open subset of $W$. Hence $D$ 
must be a tangent distribution on $W$, and Proposition \ref{nowhere} 
gives that $f|_W$ is $d$-ruled on connected components of an open 
dense subset of $W$. Moreover, the dimension of the rulings is 
bounded from below by $n+1-\dim\Ima(\varphi_Y)\geq n-2p+3$.
\vspace{2ex}\qed

\noindent \emph{Proof of Theorem \ref{local2}:} We work on the open 
dense subset of $M^n$ where $f$ is \mbox{$1$-regular} on any connected
component. Consider an open subset of a connected component where the 
index of relative nullity is $\nu\leq n-2p-1$ at any point. Lemma 
\ref{main} applies and thus the flat bilinear form $\hat\theta$ in 
\eqref{deftheta1} decomposes at any point as 
$\hat\theta=\theta_1+\theta_2$ 
where $\theta_1$ is as in part $(i)$ of that result. Hence, on any open 
subset where the dimension of $\Sal(\theta_1)
=\Sal(\hat\theta)\cap\Sal(\hat\theta)^\perp$ is constant  there are 
smooth local unit vector fields $\zeta_1,\zeta_2\in N_1$ such that 
$(\zeta_1,\zeta_2)\in \Sal(\theta_1)$. Equivalently,
\be\label{almost}
\<\beta(X,Y),\zeta_1+\zeta_2\>+\<\alpha(X,Y),\zeta_1-\zeta_2\>=0
\ee
for any $X,Y\in\mathfrak{X}(M)$. Then $\zeta_1+\zeta_2\neq 0$ 
since otherwise $\zeta_1-\zeta_2\in N_1^\perp$. Hence $\tau$ 
satisfies the condition $(*)$ and the proof follows from 
Corollary \ref{localgen}.\qed

\section{The global result}

The first two results are of independent interest.

\begin{proposition} 
Let $\tau$ be an infinitesimal bending of 
$f\colon M^n\to\R^{n+p}$ and let  $\theta$ be the flat bilinear 
form defined by \eqref{theta}. Denote 
$\nu^*(x)=\dim\Delta^*(x)$ at $x\in M^n$ where 
$$
\Delta^*(x)=\mathcal{N}(\theta)(x)=\Delta\cap\mathcal{N}(\beta)(x).
$$
Then, on any open subset of $M^n$ where $\nu^*$ is constant 
the distribution $\Delta^*$ is totally geodesic and its leaves are 
mapped by $f$ onto open subsets of affine subspaces of $\R^{n+p}$.
\end{proposition}

\proof From \eqref{parte} we have $\Delta\subset\mathcal{N}(\Y)$.  
Then \eqref{casicodazzi} and the Gauss equation give
$$
(\nap_X\beta)(Z,Y)=(\nap_Z\beta)(X,Y)=0
$$
for any $X,Y\in\Gamma(\Delta^*)$ and $Z\in\mathfrak{X}(M)$. 
Let $\nabla^*=(\nap,\nap)$ be the compatible connection in 
$N_fM\oplus N_fM$. 
Hence
$$
 0=(\nabla^*_X\theta)(Z,Y)= \theta(Z,\nabla_X Y)
$$ 
for any $X,Y\in\Gamma(\Delta^*)$ and $Z\in\mathfrak{X}(M)$.
Thus $\Delta^*\subset\Delta$ is totally geodesic.\vspace{2ex}\qed

On an open subset of $M^n$ where $\nu^*>0$ is constant consider 
the orthogonal splitting $TM=\Delta^*\oplus E$ and the 
tensor $C\colon \Gamma(\Delta^*)\times\Gamma(E)\to\Gamma(E)$
defined by
$$
C(S,X)=C_S X=-(\nabla_X S)_E
$$
where $S\in\Gamma(\Delta^*)$ and $X\in\Gamma(E)$. Since 
$\Delta^*\subset\Delta$ is totally geodesic,  the Gauss 
equation gives
$$
\nabla_TC_SX=C_SC_T +C_{\nabla_T S}
$$
for any $S,T\in\Gamma(\Delta^*)$. In particular, we have
\be\label{split}
\frac{D}{dt}C_{\gamma'}=C_{\gamma'}^2
\ee
along a unit speed geodesic $\gamma$ contained in a leaf 
of $\Delta^*$.
\medskip

The next result provides a way to transport information 
along geodesics contained in leaves of the nullity of $\theta$. 
This technique has been widely used, for instance, see \cite{DG}, 
\cite{FG} and \cite{Ji}. 

\begin{proposition}\label{completeness}
Let $\nu^*>0$ be constant on an open subset $U\subset M^n$. 
If $\gamma\colon[0,b]\to M^n$ is a unit speed geodesic such that 
$\gamma([0,b))$ is contained in a leaf of $\Delta^*$ in $U$, then 
$\Delta^*(\gamma(b))=\mathcal{P}_0^b(\Delta^*(\gamma(0)))$ where 
$\mathcal{P}_0^t$ is the parallel transport along  $\gamma$ from 
$\gamma(0)$ to $\gamma(t)$. In particular, we have
$\nu^*(\gamma(b))=\nu^*(\gamma(0))$ and the tensor $C_{\gamma'}$ 
extends smoothly to $[0,b]$. 
\end{proposition}

\proof We mimic the proof of Lemma $27$ in \cite{FG}. Let the tensor 
$J\colon E\to E$ be the solution in $[0,b)$ of
$$
\frac{D}{dt}J+C_{\gamma'}\circ J=0
$$
with initial condition $J(0)=I$. We have from \eqref{split} that 
$D^2J/dt^2=0$, and hence $J$ extends smoothly to $\mathcal{P}_0^b(E(0))$ 
in $\gamma(b)$.
Let $Y$ and $Z$ be parallel vector fields along $\gamma$ such that 
$Y(t)\in E(t)$ for each $t\in[0,b)$. Since $\gamma'\in\Delta^*$, 
it follows from \eqref{casicodazzi} that
$$
(\nabla^*_{\gamma'}\theta)(JY,Z)=(\nabla^*_{JY}\theta)(\gamma',Z).
$$
This and the definition of $J$ imply that $\theta(JY,Z)$ is parallel 
along $\gamma$. In particular $J$ is invertible in $[0,b]$. By 
continuity 
$\mathcal{P}_0^b(\Delta^*(\gamma(0)))\subset\Delta^*(\gamma(b))$,
and since $Z(0)$ is arbitrary, then 
$\mathcal{P}_0^b(\Delta^*(\gamma(0)))=\Delta^*(\gamma(b))$. 
Finally we extend the tensor $C_{\gamma'}$ to $[0,b]$
as $-DJ/dt\circ J^{-1}$.\vspace{2ex}\qed

\begin{lemma}\label{compact}
Let $f\colon M^n\to\R^{n+p}$, $p\leq 5$ and $n>2p$ be an 
isometric immersion of a compact Riemannian manifold and 
let $\tau$ be an infinitesimal bending of $f$. Then, at 
any $x\in M^n$ there is a pair of vectors 
$\zeta_1,\zeta_2\in N_fM(x)$ of unit length such that 
$(\zeta_1,\zeta_2)\in(\mathcal{S}(\theta))^\perp(x)$ where
$$
\mathcal{S}(\theta)(x)=\spa\,\{\theta(X,Y):X,Y\in T_xM\}.
$$
Moreover, on any connected component of an open dense subset of
$M^n$ the pair $\zeta_1,\zeta_2$ at $x\in M^n$ extend to smooth 
vector fields $\zeta_1$ and $\zeta_2$ parallel along $\Delta^*$ 
that satisfy the same conditions.
\end{lemma}

\proof We claim that the subset of points $U$ of $M^n$ where  
there is no such a pair, that is, where the metric induced on 
$(\mathcal{S}(\theta))^\perp$ is positive or negative definite, 
is empty. It is not difficult to see that $U$ is open. From 
Lemma~\ref{main} we have $\nu^*>0$ in $U$. 
Let $V\subset U$ be the open subset where $\nu^*=\nu_0^*$ is minimal. 
Take $x_0\in V$ and a unit speed geodesic $\gamma$ in $M^n$ contained 
in a maximal leaf of $\Delta^*$ with $\gamma(0)=x_0$. Since $M^n$ is 
compact, there is $b>0$ such that $\gamma([0,b))\subset V$ and 
$\gamma(b)\notin V$. Proposition \ref{completeness} gives 
$\nu^*(\gamma(b))=\nu_0^*$ which implies $\gamma(b)\notin U$. 
Hence, there are unit vectors $\zeta_1,\zeta_2\in N_fM(\gamma(b))$ 
such that 
$(\zeta_1,\zeta_2)\in(\mathcal{S}(\theta))^\perp(\gamma(b))$.

Let $\zeta_i(t)$ be the parallel transport along $\gamma$ of $\zeta_i$, 
$i=1,2$. Then
$$
\<\!\<\theta(X,Y),(\zeta_1,\zeta_2)\>\!\>
=\<(A_{\zeta_1-\zeta_2}+B_{\zeta_1+\zeta_2})X,Y\>.
$$
It follows from \eqref{parte} and \eqref{casicodazzi} that
\be\label{der}
(\nabla_T^*\theta)(X,Y)=(\nabla_X^*\theta)(T,Y)
\ee
where $T\in\Gamma(\Delta^*)$ extends $\gamma'$ and 
$X,Y\in\mathcal{X}(M)$. Along $\gamma$ this gives
$$
\frac{D}{dt}\mathcal{C}_{\zeta_1,\zeta_2}
=\mathcal{C}_{\zeta_1,\zeta_2}C_{\gamma'}
=C'_{\gamma'}\mathcal{C}_{\zeta_1,\zeta_2}
$$
where $\mathcal{C}_{\zeta_1,\zeta_2}
=A_{\zeta_1-\zeta_2}+B_{\zeta_1+\zeta_2}$
and  $C'_{\gamma'}$ denotes the transpose of $C_{\gamma'}$.
Moreover, by Proposition \ref{completeness} this ODE holds on $[0,b]$.
Given that $\mathcal{C}_{\zeta_1,\zeta_2}(\gamma(b))=0$, then 
$\mathcal{C}_{\zeta_1,\zeta_2}$ vanishes along 
$\gamma$. This is a contradiction and proves the claim.

We have from \eqref{der}  that
$$
(\nabla_T^*\theta)(X,Y)
=-\theta(\nabla_XT,Y)\in\Gamma(\mathcal{S}(\theta))
$$
for any $T\in\Gamma(\Delta^*)$ and $X,Y\in\mathcal{X}(M)$. Thus
$\mathcal{S}(\theta)$ is parallel along the leafs of $\Delta^*$.
Let $U_0$ be a connected component of the open dense subset of $M^n$
where the dimension of $\Delta^*$, $\mathcal{S}(\theta),
\mathcal{S}(\theta)\cap \mathcal{S}(\theta)^\perp$ and the index of
the metric induced on 
$\mathcal{S}(\theta)^\perp\times\mathcal{S}(\theta)^\perp$ are all
constant. Hence on $U_0$ the vector fields $\zeta_1,\zeta_2$ can be 
taken parallel along the leafs of $\Delta^*$. \vspace{2ex}\qed
  
For an hypersurface  $f\colon M^n\to\R^{n+1}$ we have
\be\label{defB1}
(\nab_XL)Y=\<B_NX,Y\>N+f_*\mathcal{Y}(X,Y)
\ee
where $N$ is a unit vector field normal to $f$. The next result follows 
from Theorem $13$ in \cite{DV} and was fundamental in \cite{Ji}.

\begin{lemma}\label{triv} An infinitesimal bending $\tau$ of 
$f\colon M^n\to \R^{n+1}$ is trivial if and only if $B_N=0$.
\end{lemma} 

\noindent \emph{Proof of Theorem \ref{maincompact}:} 
We assume that there is no open subset of $M^n$ where the index 
of relative nullity satisfies $\nu\geq n-1$. By Lemma \ref{compact}, 
on connected components of an open dense subset of $M^n$ there are 
$\zeta_1,\zeta_2\in\Gamma(N_fM)$ with 
$\|\zeta_1\|=\|\zeta_2\|=1$ parallel along the leaves of 
$\Delta^*$ and such that 
$$
\<\!\<\theta(X,Y),(\zeta_1,\zeta_2)\>\!\>=0
$$
for any $X,Y\in \mathfrak{X}(M)$. It follows from \eqref{theta} 
that \eqref{almost} holds on connected components of an open dense 
subset of $M^n$.
Let $U\subset M^n$ be an open subset where $\zeta_1, \zeta_2$ are
smooth and $\zeta_1+\zeta_2\neq 0$. Thus $\tau|_U$ satisfies the 
condition~$(*)$. Let $\tilde{V}\subset U$ be an open subset where 
$\tau$ is a genuine infinitesimal bending.
By Corollary~\ref{localgen} we have that $f$ is $(n-1)$-ruled on 
each connected component $V$ of an open dense subset of $\tilde{V}$. 
Since our goal is to show that $V$ is empty we assume 
otherwise.

Proposition \ref{nowhere} and the proof of Theorem \ref{local} 
yield that the rulings on $V$ are determined by the tangent subbundle 
$D=\ker\varphi_Y$ where $\varphi$ was given in Lemma \ref{varphi}
and $Y\in RE(\varphi)$. Also from that proof
$\dim\Ima(\varphi_Y)=2$ and therefore
$\Ima(\varphi_Y)=R\oplus R$ where $N_fM=P\oplus R$ as in 
Lemma \ref{varphi}. Lemma~\ref{fbn} gives
$$
\varphi_X(D)\subset\Ima(\varphi_Y)\cap\Ima(\varphi_Y)^\perp=\{0\}
$$
for any $X\in\mathfrak{X}(M)$, that is, $D=\mathcal{N}(\varphi)$. 
In particular, from the definition of $\varphi$ it follows that 
$D\subset \mathcal{N}(\alpha_R)$. Hence, by dimension reasons either 
$\mathcal{N}(\alpha_R)=TM$ or $D=\mathcal{N}(\alpha_R)$.
Next we contemplate both possibilities.

Let $V_1\subset V$ be an open subset where $\mathcal{N}(\alpha_R)=TM$
holds, that is, $N_1=P$.  Thus $N_1$ is parallel 
relative to the normal connection since, otherwise, the Codazzi 
equation gives $\nu=n-1$, and that has been ruled out.
Hence $f|_{V_1}$ reduces codimension, that is, $f(V_1)$ is contained 
in an affine hyperplane $\R^{n+1}$.
Decompose $\tau=\tau_1+\tau_2$ where $\tau_1$ and $\tau_2$ 
are tangent and normal to $\R^{n+1}$, respectively. It follows that 
$\tau_1$ is an infinitesimal bending of $f|_{V_1}$ in $\R^{n+1}$. 
Since $\tau$ satisfies the condition  $(*)$ then  
Lemma \ref{triv} gives that $\tau_1$ is trivial, that is, the restriction 
of a Killing vector field of $\R^{n+1}$ to $f(V_1)$. Extending $\tau_2$ as 
a vector field normal to $\R^{n+1}$ it follows that $\tau|_{V_1}$ extends 
in the singular sense and this is a contradiction.

Let $V_2\subset V$ be an open subset where $D=\mathcal{N}(\alpha_R)$. 
By assumption $D\neq \Delta$. Let $\hat{D}$ be the distribution 
tangent to the rulings in a neighborhood $V_2'$ of $x_0\in V_2$. 
From Proposition \ref{nowhere} we have  $D(x_0)=\hat{D}(x_0)$. 
Let $W\subset V_2'$ be an open subset where $D\neq\hat{D}$, that is, 
where $D$ is not totally geodesic. Then there are two transversal 
$(n-1)$-dimensional rulings passing through any point $y\in W$. 
It follows easily that $N_1=P$ on $W$. As above we obtain that 
$\tau|_W$ extends in the singular sense, leading to a contradiction. 
Let $V_3\subset V_2$ be the interior of the subset where $D$ is 
totally geodesic. On $V_3$ the Codazzi equation gives 
$$
\nabla^\perp_X\alpha(Z,Y)\in \Gamma(P)
$$
for all $X,Y\in \Gamma(D)$ and $Z\in\mathcal{X}(M)$. Thus $R$ 
is parallel along $D$ relative to the normal connection. 
We have from Proposition 4 in \cite{DG}  that $f$ admits a 
singular extension 
$$
F(x,t)=f(x)+t\lambda(x)
$$ 
for $\lambda\in\Gamma(f_*TM\oplus P)$ as a flat hypersurface. 
Moreover, $F$ has $R$ as normal bundle and $\partial_t$ belongs 
to the relative nullity distribution. 
Then $(\nab_X\lambda)_R=0$ for any $X\in\mathfrak{X}(V_3)$. Hence 
\eqref{requisito} is satisfied and thus $\tau|_{V_3}$ extends in 
the singular sense. This is a contradiction which shows that $V$ 
is empty, and hence also is $\tilde{V}$.

It remains to consider the existence of an open subset 
$U'\subset M^n$ where $\zeta_1,\zeta_2$ are smooth and 
$\zeta_1+\zeta_2=0$. It follows from \eqref{almost} that 
$\zeta_1-\zeta_2\perp N_1$. Once more, we obtain that 
$f(U')\subset\R^{n+1}$. Thus, we have an orthogonal
decomposition of $\tau|_{U'}$ as in part $(ii)$ of the 
statement and $\tau_1,\tau_2$ extend in the singular sense 
as follows:
\begin{itemize}
\item[(i)] 
$\bar{\tau}_1(x,t)=\tau_1(x)$ to $F\colon U\times\R\to\R^{n+2}$ 
where $F(x,t)=f(x)+te$. 
\item[(ii)] For instance locally as $\bar{\tau}_2(x,t)=\tau_2(x)$ 
to $F\colon U\times I\to\R^{n+2}$ where $F(x,t)=f(x)+tN$
being $N$ is a unit normal field to $f|_U$ in $\R^{n+1}$.\qed
\end{itemize}

\begin{remarks}{\em $(1)$ In case $(ii)$ of Theorem \ref{maincompact} 
if $\tau_1$ is trivial then $\tau_1$ and $\tau_2$ extend in the same 
direction, and hence $\tau$ also does. Therefore we are also in case 
$(i)$.

\noindent $(2)$ Notice that for $p=2$ we have shown as part of the 
proof that an infinitesimal bending of a submanifold without flat 
points as in in part $(ii)$ of Theorem \ref{local2} cannot be genuine.
}\end{remarks}

\section{Nonflat ambient spaces}

In this section we argue for the following statement: 
\medskip

\noindent \emph{Theorems \ref{local1}, \ref{local2} 
and \ref{local} hold if the Euclidean ambient space is replaced 
by a nonflat space form.}
\vspace{1ex}

Let $f\colon M^n\to\Q^{n+p}_c$ be an isometric immersion where 
$\Q^{n+p}_c$ denotes either the sphere  $\Sf_c^{n+p}$ or the 
hyperbolic space $\Hy_c^{n+p}$ of sectional curvature $c\neq 0$.
Then we say that $\tau\in\Gamma(f^*T\Q_c^{n+p})$ is an 
infinitesimal bending of $f$ if \eqref{infbend} is satisfied 
in terms of the connection in $\Q_c^{n+p}$. And now that $f$ is 
\emph{$r$-ruled} means that there is an $r$-dimensional smooth 
totally geodesic distribution whose leaves are mapped by $f$ to 
open subsets of totally geodesic submanifolds of the 
ambient space $\Q^{n+p}_c$.
\medskip

In the sequel, for simplicity we also denote by $f$  the composition 
of the immersion with the umbilical inclusion of $\Q_c^{n+p}$ into 
$\mathbb{O}^{n+p+1}$, where $\mathbb{O}^{n+p+1}$ stands for either 
Euclidean or Lorentzian flat space depending on whether $c>0$ or 
$c<0$, respectively.  

Let $\tau$ be an infinitesimal bending of $f$ and let 
$\F\colon I\times M^n\to\Q_c^{n+p}$ be a smooth variation such 
that $f_t=\F(t,\cdot)\colon M^n\to\Q_c^{n+p}$ verifies $f_0=f$ 
and having $\tau$ as variational vector field. In this case
we still have that \eqref{der metrica}, \eqref{der conex} and 
\eqref{der curv} hold. And also as before, associated to $\tau$ 
we have the tensors 
$$
LX=\nab_{X}\tau\;\;\mbox{and}\;\;B(X,Y)=(\nab_XL)Y
$$
where $X,Y\in\mathfrak{X}(M)$ and $\nab$ denotes the connection 
in $\Q_c^{n+p}$. Now
$$
B(X,Y)=f_*\Y(X,Y)+\beta(X,Y)+c\<f_*Y,\tau\>f_*X-c\<X,Y\>\tau
$$
where the tensors $\Y\colon TM\times TM\to TM$ and 
$\beta\colon TM\times TM\to N_fM$ are the tangent and normal 
components of $\partial/\partial t|_{t=0}\alpha^t$, respectively, 
and $\alpha^t$ is the second fundamental form of $f_t$ as 
a submanifold in $\Q^{n+p}_c$. In particular, we have that
\eqref{derGauss} holds.

In this case, an infinitesimal bending of $f$ is said to satisfy the 
\emph{condition $(*)$} if there is $\eta\in\Gamma(N_fM)$ of unit 
length and $\xi\in\Gamma(R)$, where $R$ is determined by the orthogonal 
splitting $N_fM=P\oplus R$ and $P=\spa\{\eta\}$, such that
$$
B_\eta+A_\xi+c\<\tau,\eta\>\id=0
$$
where $B_\eta=\<\beta,\eta\>$.  

The cone over an isometric immersion $f\colon M^n\to\Q^{n+p}_c$ is 
defined by
\begin{align*}
\hat{f}\colon\hat{M}^{n+1}
=(0,\infty)\times &M^n\to\mathbb{O}^{n+p+1}\nonumber\\
(s,x)&\mapsto sf(x).
\end{align*}
Notice that $\partial_s$ lies in the relative nullity of $\hat{f}$ 
and that $N_{\hat{f}}\hat{M}$ is the parallel transport of $N_fM$ 
along the lines parametrized by $s$. Observe that if $c<0$, then
the cone over $f$ is a Lorentzian submanifold of $\mathbb{L}^{n+p+1}$ 
and hence $N_{\hat f}\hat{M}$ has positive definite metric.

If $\tau$ is an infinitesimal bending of $f$, it is easy to see
that $\hat{\tau}(s,x)=s\tau(x)$ is an infinitesimal 
bending of $\hat{f}$ in $\mathbb{O}^{n+p+1}$, that is, 
$\hat{\tau}$ is a vector field that satisfies \eqref{infbend} with 
respect to the connection in $\mathbb{O}^{n+p+1}$.
Moreover, if $\tau$ satisfies the condition $(*)$ then  $\hat{\tau}$ 
satisfies the condition $(*)$ for the flat ambient space.

Let $\hat{f}$ be the cone over an immersion $f$ in $\Q_c^{n+p}$. 
Notice that the parameter $s$ defines lines parallel to the position 
vector. Thus, if the map $\hat{f}+t\lambda$, is a singular extension of 
$\hat{f}$ for some vector field $\lambda$ then the intersection 
of its image with $\Q_c^{n+p}$ determines a singular extension of $f$.

Consider the maps  
$$
\hat{F}(t,s,x)=\hat{f}(s,x)+t\lambda(s,x)\;\;\mbox{and}\;\;
\hat{\tau}'(t,s,x)=\hat{\tau}(s,x)+t\bar{L}\lambda(s,x)
$$
as in the proofs of Theorems \ref{local1} and \ref{local}. 
Notice that
\begin{align*}
\<\hat{F}(t,s,x),\hat{\tau}'(t,s,x)\>
&=\<\hat{f}(s,x)+t\lambda(s,x),\hat{\tau}(s,x)
+t\bar{L}\lambda(s,x)\> \\
&=st\<f(x),\bar{L}\lambda\>+st\<\lambda,\tau\>\\
&=0
\end{align*}
where for the last equality we used $\hat{L}\partial_s=\tau(x)$. 
Then we have that $\hat{\tau}'$ is orthogonal to the position 
vector $\hat{F}$. From this we have that if $\hat{F}$ determines 
a singular extension of $\hat{f}$ then $\tau$ extends in the 
singular sense. 

As in the proofs of Theorems \ref{local1} and \ref{local}, if 
there is no $\lambda$ as above that determines a singular extension 
of $\hat{f}$ we conclude that $\hat{f}$ is ruled. Finally, observe 
that being $\hat{f}$ the cone over $f$, then these rulings determine 
rulings of $f$.

\noindent Marcos Dajczer\\
IMPA -- Estrada Dona Castorina, 110\\
22460--320, Rio de Janeiro -- Brazil\\
e-mail: marcos@impa.br

\bigskip

\noindent Miguel Ibieta Jimenez\\
IMPA -- Estrada Dona Castorina, 110\\
22460--320, Rio de Janeiro -- Brazil\\
e-mail: mibieta@impa.br
\end{document}